\documentclass[1p,10pt]{elsarticle}

\journal{arXiv.org}

\usepackage{amssymb}
\usepackage{amsmath}
\usepackage{graphicx}
\usepackage[T1]{fontenc}

\begin{document}

\newcommand{\BR}{{\mathbb R}}
\newcommand{\BN}{{\mathbb N}}
\newcommand{\BZ}{{\mathbb Z}}
\newcommand{\BC}{{\mathbb C}}
\newcommand{\e}{{\bf e}}
\renewcommand{\a}{{\bf a}}
\renewcommand{\b}{{\bf b}}
\renewcommand{\u}{{\bf u}}

\newcommand{\m}{{\bf m}}
\newcommand{\p}{{\bf p}}

\newcommand{\f}{{\bf f}}
\newcommand{\g}{{\bf g}}

\newcommand{\w}{{\bf w}}
\newcommand{\G}{{\bf G}}
\newcommand{\K}{{\bf K}}

\newcommand{\cl}{C \kern -0.1em \ell}

\newcommand{\ov}{\overline}
\newcommand{\un}{\underline}
\newcommand{\proof}{\bf {Proof:} \rm}
\renewcommand{\qed}{$\blacksquare$}
\newtheorem{theorem}{Theorem}[section]
\newtheorem{remark}{Remark}[section]
\newtheorem{lemma}{Lemma}[section]
\newtheorem{proposition}{Proposition}[section]
\newtheorem{corollary}{Corollary}[section]
\newtheorem{definition}{Definition}[section]
\newtheorem{example}{Example}[section]

\begin{frontmatter}
 \title{Classes of hypercomplex polynomials of discrete variable based on
the quasi-monomiality principle}
\author{N.~Faustino\fnref{label2}\corref{cor1}}
\ead{faustino@ime.unicamp.br}
\ead[url]{https://sites.google.com/site/nelsonfaustinopt/}
\cortext[cor1]{Corresponding author}
\fntext[label2]{N. Faustino was supported by fellowship 13/07590-8 of FAPESP (S.P., Brazil).}
\address{Departamento de Matem\'atica Aplicada, IMECC--Unicamp,
CEP 13083--859, Campinas, SP, Brasil}


\begin{abstract}
 With the aim of derive a quasi-monomiality formulation in the context of discrete hypercomplex variables, 
 one will amalgamate through a Clifford-algebraic structure of signature $(0,n)$ the umbral calculus framework
 with Lie-algebraic symmetries. The exponential generating function ({\bf EGF}) carrying the {\it continuum}
 Dirac operator $D=\sum_{j=1}^n\e_j\partial_{x_j}$ 
 together with the Lie-algebraic representation of raising and lowering operators acting on the lattice $h\BZ^n$ 
 is used to derive the corresponding hypercomplex polynomials of discrete variable as Appell sets
 with membership
 on the space Clifford-vector-valued polynomials. Some particular
 examples concerning this construction such as the hypercomplex versions of falling factorials
 and the Poisson-Charlier polynomials
 are introduced.
 Certain applications from the view of interpolation theory and integral transforms are also discussed. 
\end{abstract}

\end{frontmatter}


\section{Introduction}

\subsection{State of Art}

The modern language of Appell/Sheffer sets, described in terms of the {\it
quasi-monomiality} formalism (see e.g.
\cite{DattoliSriKhan05},~\cite{BelDattoliKhanRicci07} and the
references given there), plays nowadays a central role in the
construction of polynomial solutions for finite difference
equations possessing Lie-algebraic symmetries (cf.
\cite{SmirnovTurbiner95},\cite{DHS96} \& \cite{LTW04}).
Although this formalism has been fully developed and popularized by
Ben Cheikh, Dattoli, Srivastava {\it et all} (see e.g.
\cite{BenCheikhZagh03},~\cite{SriBenCheik03},~\cite{DattoliSriKhan05},
\cite{BlasDattoliHorzPenson06} and the references given there), the fundamentals of such theory were
sketched in the former paper of Di Bucchianico-Loeb-Rota \cite{BLR98}.

Indeed, the construction of polynomial solutions based on operational methods has a long history
that started with the study of classical equations of motion from
the quantum mechanical side by Wigner (cf. \cite{Wigner50}). Years later,
Turbiner-Ushveridze \cite{TurbinerUsh87}, Gagnon-Winternitz
\cite{GagnonWinternitz88} and numerous followers across extend this characterization,
the so-called {\it quasi-exact solvability} condition, to generate polynomial spaces 
as finite-dimensional irreducible Fock spaces that invariant under the action of a {\it degree-preserving} operator,
In the setting of quantum mechanics such operator may be described as a number operator
(cf.~\cite{BlasDattoliHorzPenson06},\cite{VinetZhedanov09}). 

In the context of hypercomplex variables,
this quantum mechanical correspondence provided by the {\it quasi-exact solvability} condition
was successfully applied in \cite{CFK11} to describe the spectrum of the
quantum harmonic oscillator in terms of Clifford-Hermite
functions/polynomials (see, for instance, \cite[Theorem 3.1 \& Theorem 3.2]{CFK11}).

The main goal of this paper is the hypercomplex extension of the
quasi-monomiality formalism to discrete variables in
interrelationship with Lie-algebraic representations of finite
difference operators which are Clifford-vector-valued. In {\it
continuum} this is currently an ongoing research topic (cf.
\cite{CFM11},~\cite{Eelbode12},\cite{Lavicka12}).

\subsection{The Lie-algebraic background through hypercomplex variables}\label{LieAlgebraicSubsection}

 As described in \cite{FR11}, the main idea around Lie-algebraic discretization in the context of
 hypercomplex variables consists into start from a given algebra of radial type acting on the
 space of Clifford-vector-valued
 polynomials $\mathcal{P}=\BR[x]\otimes
\cl_{0,n}$, a Hilbert module generated from the tensor product between the ring 
of multi-variable polynomials $\BR[x]$ ($x\in \BR^n$) and the Clifford algebra
$\cl_{0,n}$ of signature $(0,n)$. 

The rigorous formulation of this approach is based on the study of left representations
with membership on the algebra of
endomorphisms $\mbox{End}(\mathcal{P})$ through the
canonical correspondence
\begin{eqnarray}\label{RadialAlgebra}
\mbox{End}(\mathcal{P})\cong \mbox{Alg}\left\{
L_j,M_j,\e_j~:~j=1,\ldots,n\right\}.
\end{eqnarray}

Here and elsewhere $L_1,L_2,\ldots,L_n,M_1,M_2,\ldots,M_n$ together with the identity operator
$I:\f(x) \mapsto \f(x)$ are assumed to be the canonical generators of the Weyl-Heisenberg algebra with
dimension $2n+1$ satisfying the graded commuting relations
\begin{eqnarray}
\label{WeylHeisenbergLM}
\begin{array}{llll}
[L_j,L_k]=0,& \left[M_j,M_k\right]=0, &
\left[L_j,M_k\right]=\delta_{jk}I
\end{array}
\end{eqnarray}
whereas $\e_1,\e_2,\ldots,\e_n$ is fixed as an orthogonal basis of $\BR^n$
satisfying, for each $j,k=1,2,\ldots,n$, the set of graded anti-commuting relations
\begin{eqnarray}
\label{CliffordGenerators} \e_j\e_k+\e_k\e_j=-2\delta_{jk}.
\end{eqnarray}

In this context, the canonical description of $\mbox{End}(\mathcal{P})$
provided by (\ref{RadialAlgebra}) thus corresponds to the Lie-algebraic counterpart
of the radial algebra introduced by Sommen in \cite{Sommen97}. Moreover, the left endomorphisms 
$L=\sum_{j=1}^n \e_j L_j$ and $M=\sum_{j=1}^n \e_j
M_j$ encode the underlying symmetries of the orthosymplectic Lie algebra
$\mathfrak{osp}(1|2)$ (cf. \cite[Subsection 2.3]{FR11}).

As examples of Appell sequences that can be derived from this scheme one can mention
the hypercomplex extension of the multi-index factorial powers considered in
\cite[Subsection 3.2]{FK07} as the discrete counterpart of the
Clifford-vector-valued homogeneous polynomials used to determine the
Fischer decomposition with respect to the finite difference Dirac
operators $D_h^\pm$ of forward/backward type (see also \cite[Subsection 3.1]{FR11}). In intertwine with
Howe's invariant theory \cite{Howe89TransAMS},
the same scheme was recently exploited by the author in \cite{Faustino13} to construct {\it degree-preserving}
spaces of Clifford-vector-valued
polynomials on the lattice $h\BZ^n$ as invariant and irreducible spaces with respect to the
Howe dual pair $(SO(n),\mathfrak{su}(1,1))$.

The {\it quasi-monomiality} principle\footnote{From the view of quantum field theory, the {\it quasi-monomiality}
principle provides a meaningful interpretation of the second quantization formalism.} provided by
the set of operators
$L_j$ and $M_j$ amalgamates
most of the substantial work already done in the study of multi-variable Appell sequences
through continuous and discrete variables
(cf.~\cite{BenCheikhZagh03},\cite{DattoliSriKhan05},\cite{BelDattoliKhanRicci07},\cite{VinetZhedanov09}). 
In concrete, if the set of multi-variable polynomials
$\{\m_\alpha(x)~:~\alpha \in \mathbb{N}_0^n\}$, determined by the
condition $\m_{{\bf 0}}(x)=1$ (${\bf 0}\in \BR^n$) and by the set of quasi-monomiality constraints
\begin{eqnarray} \label{quasiMonomialOp}
L_j\m_{\alpha}(x)=\alpha_j\m_{\alpha-\e_j}(x)& \mbox{and}&
M_j\m_{\alpha}(x)=\m_{\alpha+\e_j}(x),
\end{eqnarray}
with $\alpha=\sum_{j=1}^n \alpha_j\e_j$ has
\begin{eqnarray}
\label{EGFSheffer} \sum_{|\alpha|=0}^\infty
\m_\alpha(x)~\frac{y^\alpha}{\alpha!} &=&\prod_{j=1}^n
  \dfrac{1}{\kappa\left(\ell^{-1}(y_j)\right)}~\exp\left(x_j\ell^{-1}(y_j)\right)
\end{eqnarray}
as multi-variable exponential generating function ({\bf EGF})\footnote{In the
operational form, the multi-variable {\bf EGF}~may be represented through the action
of the exponentiation operator $\exp\left(\sum_{j=1}^n y_j
M_j\right)$ on $\m_{{\bf 0}}(x)=1$.} then, the associated set of ladder
operators $L_j$ and $M_j$, given by (cf. \cite[Section
3]{BlasDattoliHorzPenson06})
\begin{eqnarray}
\label{LMSheffer}
\begin{array}{lll}
L_j&=&\ell\left(\partial_{x_j}\right) \\
M_j&=&\left(x_j-\kappa'\left(\partial_{x_j}\right)\kappa\left(\partial_{x_j}\right)^{-1}\right)
\ell'\left(\partial_{x_j}\right)^{-1}
\end{array}
\end{eqnarray}
satisfy the Weyl-Heisenberg graded commuting relations (\ref{WeylHeisenbergLM}).

The implicit mathematical conditions corresponding to the above
characterization (see also \cite[Theorem 2.5.3]{Roman84} and
\cite[Theorem 3.7.1]{Roman84}) are the isomorphism between the
algebra of formal power series and the algebra of linear
functionals carrying the algebra of multi-variable polynomials
(cf. \cite[Theorem 2.1.1]{Roman84}) and the {\it shift-invariant} property (cf.
\cite[Corollary 2.2.8]{Roman84}) underlying the set of operators $\kappa\left(\partial_{x_j}\right)$ and
$\ell\left(\partial_{x_j}\right)$ obtained through the substitutions
$t \rightarrow \partial_{x_j}$ on
\begin{eqnarray*}
\kappa(t)=\sum_{k=0}^\infty a_k \frac{t^k}{k!}& \mbox{and} &
\ell(t)=\sum_{k=1}^\infty b_k \dfrac{t^k}{k!}.
\end{eqnarray*}

The one-to-one correspondence between linear functionals and formal power series also shows that
$\ell'\left(\partial_{x_j}\right)$ is a {\it shift-invariant} operator that
coincides with $L_j':=[L_j,x_j]$, the so-called {\it
Pincherle derivative}\footnote{In Roman's book \cite{Roman84} the
Pincherle derivative associated to
$\ell'\left(\partial_{x_j}\right)$ can be found on {\bf Chapter 2}, {\bf
Section 6} under the name of {\it umbral
shift}.} (cf. \cite[Section 2]{DHS96}). The
existence of $\ell'\left(\partial_{x_j}\right)^{-1}$ is thus
assured by the condition $\ell'\left(\partial_{x_j}\right)1 \neq
0$ (cf. \cite[Subsection 2.1]{LTW04}).

\subsection{Outline of the paper}

In this paper, the {\it quasi-monomialy} approach is extended to discrete hypercomplex variables
based on the Lie-algebraic formulation revisited in Subsection \ref{LieAlgebraicSubsection}.
To begin with, we started Section \ref{QuasiMonomialitySection} by formulating
the quasi-monomiality principle with respect to the finite difference Dirac operator $D_h^+$
and the multi-variable {\bf EGF} of the form
\begin{eqnarray}
\label{generalEGFhypergeometric} G_h(x,y;\kappa)=\prod_{j=1}^n
\dfrac{1}{\kappa\left(
\frac{1}{h}\log\left( 1+h y_j\right)\right)}~\left( 1+h y_j\right)^{\frac{x_j}{h}}.
\end{eqnarray}

Based on the knowledge of (\ref{generalEGFhypergeometric}), the main results concerning this approach
are also proved in Section \ref{QuasiMonomialitySection}. One starts in Proposition \ref{FourierDualProp}
by showing that, under the existence of a multi-variable function $\lambda(y)$ ($y\in \BR^n$)
satisfying the constraint 
$$\lambda\left( \frac{D_h^+\exp(x\cdot y)}{\exp(x\cdot y)}\right)=\prod_{j=1}^n\kappa(y_j),$$
the Fourier dual $\Lambda_h$ of $D_h^+$ is uniquely determined. 
Afterwards, based on the knowledge of $\lambda(y)$ it follows from Proposition \ref{RodriguesFormulaXh}
that $\Lambda_h$
admits an operational representation similar to the Rodrigues formula used in \cite[Theorem 3.2]{CFK11} 
to determine an operational representation for the Clifford-Hermite polynomials. As briefly sketched in 
Remark \ref{HypercomplexPoissonCharlierRemark}, the hypercomplex version of the Poisson-Charlier polynomials
obtained in Example \ref{HypercomplexPoissonCharlier}
may be viewed as discrete analogues for the Clifford-Hermite polynomials.

Proposition \ref{RodriguesFormulaXh} is also applied to show that, for a given set
of $\mu_k$'s, the resulting Appell set $\{ \w_k(x;h;\lambda)~:~k\in
\BN_0\}$ of Clifford-vector-valued 
polynomials, constructed from the operational rule 
$\w_k(x;h;\lambda)=\mu_k\left(
\Lambda_h\right)^k\a$ ($\a\in \cl_{0,n}$), may be determined as an integral type transform written 
in terms of the hypercomplex version for the multi-index falling factorials obtained in
Example \ref{HypercomplexFactorialPolynomials}.
This characterization corresponds to Corollary \ref{RodriguesFormulaXhCorollary}.

From the combination of the aformentioned results with the hypercomplex {\bf EGF} $\G(x,t)$
obtained in Proposition \ref{EGFDiracProposition} for the {\it continuum} Dirac operator
$D=\sum_{j=1}^n\e_j \partial_{x_j}$,
one gets in Corollary 
\ref{hypercomplexEGFCorollary} a closed formula for the hypercomplex {\bf EGF} $\G_h(x,t;\lambda)$ 
underlying the Appell set $\{ \w_k(x;h;\lambda)~:~k\in
\BN_0\}$ 
as the hypercomplex counterpart of the multi-variable {\bf EGF} (\ref{generalEGFhypergeometric}).

As an illustration of this approach, we discuss in Section \ref{ConcludingSection}
further applications in the context of interpolation theory and integral transforms as well. 
In the same train of though of 
\cite{GuerlbeckHommel03}, 
the restriction of the continuous Fourier transform to the $n-$dimensional cube
$Q_h=\left(-\frac{\pi}{h},\frac{\pi}{h}\right)^n$
is used to get an integral representation for each Clifford-vector-valued 
polynomial $\w_k(x;h,\lambda)$ on the momentum space.

The operational framework developed throughout this paper is slighly different from
other recent approaches on this direction (see e.g. \cite[Sections 4 \& 5]{CFM11} and \cite[Section 4]{Eelbode12}),
though there is some connection with these (see Remark \ref{BesselConnection}).
The main novelty addressed here is the incorporation of the multi-variable {\bf EGF} on the {\it quasi-monomiality}
formulation. On the other hand, through the construction obtained in Example \ref{BernoulliExample}
the presented approach shall also fits into the framework developed
in \cite[Section 2]{MT08} for 'generalized' Bernoulli polynomials,
although it does not incorporates
{\it a-priori} a set of 'regular variables' (see e.g. \cite[Subsection 2.2]{GuerlebeckSproessig97} on which such
construction was revisited and discussed in depth).

\section{Quasi-Monomiality through discrete hypercomplex variables}\label{QuasiMonomialitySection}

\subsection{Discrete Clifford calculus}

In this subsection some basic definitions and properties carrying the
discrete Clifford setting will be collected. 
Further details concerning the definition and properties of Clifford algebras
may be found in \cite{Sommen97}, \cite[Chapter 1]{GuerlebeckSproessig97} or even
in \cite[Chapter 2]{RodriguesOliveira07}.
For an overview of finite difference discretizations of Dirac
operators one refer to \cite[Chapter 5]{GuerlebeckSproessig97} and \cite[Section 2]{FK07}.
For the construction of finite difference
discretizations on the lattice based on the interplay between finite difference calculus and Lie-algebraic
symmetries one refer to \cite[Subsection 2.1]{FR11} and \cite[Section 2]{Faustino13}.

Let $\BR[x]$ be the
ring of multi-variable polynomials in the variable $x\in \BR^n$ and
$\e_1,\e_2,\ldots,\e_n$ an
orthogonal basis of $\BR^n$.
The Clifford algebra of signature $(0,n)$, denoted by $\cl_{0,n}$, is an algebra with dimension $2^n$ in
which the $\e_j$'s satisfy the graded commuting relations (\ref{CliffordGenerators}).
Under the linear space isomorphic between $\cl_{0,n}$ and the
exterior algebra $\Lambda^*\left(\BR^n\right)$ determined from linearity arguments through the set of mappings
$$\e_{j_1}\e_{j_2}\ldots \e_{j_r} \mapsto dx_{j_1}dx_{j_2}\ldots dx_{j_r},$$
with $1\leq j_1<j_2<\ldots<j_r\leq n$, it follows that every Clifford number may be written
as a linear combination involving
$r$-multivector bases of the form $\e_{j_1}\e_{j_2}\ldots \e_{j_r}$.
That is, for a given subset $J=\{j_1,j_2,\ldots,j_r\}$ of $\{ 1,2,\ldots,n\}$, any $\a \in \cl_{0,n}$
is represented through the summation formula
 \begin{eqnarray*}
\a=\sum_{r=0}^n\sum_{|J|=r} a_J ~\e_J, & \mbox{with} & \e_{J}=\e_{j_1}\e_{j_2}\ldots \e_{j_r}.
\end{eqnarray*}

 In particular, any vector $(x_1,x_2,\ldots,x_n)$ of $\BR^n$ are represented as
 $x=\sum_{j=1}^n x_j \e_j$ whereas the translations $(x_1,x_2,\ldots, x_j\pm
 h,\ldots,x_n)$ on the lattice $h\BZ^n \subset \BR^n$ with mesh width $h>0$ are described
as $x\pm h\e_j$. In the same order of ideas, every multi-index $(\alpha_1,\alpha_2,\ldots,\alpha_n)$ with
membership in $\BN^n_0$ may be represented as $\alpha=\sum_{j=1}^n\alpha_j\e_j$.
 
 The Clifford algebra $\cl_{0,n}$ is indeed an associative algebra with identity $1$, 
 containing $\BR$ and $\BR^n$ as subspaces. Its center, generated from the anti-commutator $xy+yx$
 between two Clifford vectors $x=\sum_{j=1}^n x_j\e_j$ and $y=\sum_{j=1}^n y_j\e_j$, gives rise to
the following inner product relation
$$\sum_{j=1}^nx_j y_j=-\frac{1}{2}(xy+yx),$$
denoted throughout this paper by $x \cdot y$.

Next, for each Clifford-vector-valued function
$\f(x)=\sum_{r=0}^n\sum_{|J|=r} f_J(x) \e_J$, with $f_J(x)$
real-valued, one define the forward/backward discretizations for the partial
derivatives $\partial_{x_j}$ on the lattice $h\BZ^n$ as
\begin{eqnarray}
\label{DiffPmj}(\partial_{h}^{+j}\f)(x)=\dfrac{\f(x+h\e_j)-\f(x)}{h} &\mbox{and}
& (\partial_{h}^{-j}\f)(x)=\dfrac{\f(x)-\f(x-h\e_j)}{h}.
\end{eqnarray}

The forward and backward finite difference operators, $\partial_h^{+j}$ and $\partial_h^{-j}$ respectively,
are intertwined by the translation operators $(T_h^{\pm j}\f)(x)=\f(x\pm h \e_j)$ i.e.
\begin{eqnarray}
\label{TranlationsPmj}T_h^{-j}(\partial_h^{+j}\f)(x)=(\partial_h^{-j}\f)(x)
& \mbox{and} &
T_h^{+j}(\partial_h^{-j}\f)(x)=(\partial_h^{+j}\f)(x).
\end{eqnarray}

Moreover, they satisfy the set of product rules 
\begin{eqnarray}
\label{productRule}
\begin{array}{ccc}
\partial_h^{+j}\left( \g(x)\f(x)  \right)
&=&(\partial_{h}^{+ j}\g)(x)\f(x+ h\e_j)+\g(x)(\partial_h^{+
j}\f)(x)\\
\partial_h^{- j}\left( \g(x)\f(x)  \right)
&=&(\partial_{h}^{-j}\g)(x)\f(x- h\e_j)+\g(x)(\partial_h^{-
j}\f)(x).
\end{array}
\end{eqnarray}

Based on (\ref{DiffPmj}), one introduce the forward/backward
discretizations $D_h^\pm$ for the {\it continuum} Dirac operator
$D=\sum_{j=1}^n \e_j \partial_{x_j}$ as follows:
\begin{eqnarray}
\label{Dhpm} D_h^+ =\sum_{j=1}^n \e_j \partial_h^{+j} &
\mbox{and}& D_h^- =\sum_{j=1}^n \e_j \partial_h^{-j}.
\end{eqnarray}

Let us now restrict ourselves to Clifford-vector-valued functions
with membership in the space of Clifford-vector-valued polynomials $\mathcal{P}=\BR[x]\otimes
\cl_{0,n}$. From the set of product rules (\ref{productRule}) one can see
that the simplest representation underlying the algebra of endomorphisms
$\mbox{End}(\mathcal{P})$, as described through the canonical correspondence (\ref{RadialAlgebra}), is
given by the forward/backward discretizations $\partial_{h}^{+j}$ resp.
$\partial_{h}^{-j}$ and the multiplication operators
$x_jT_h^{-j}:\f(x) \mapsto x_j \f(x-h\e_j)$ resp.
$x_jT_h^{+j}:\f(x) \mapsto x_j \f(x+h\e_j)$. Indeed, a simple
computation based in (\ref{productRule}) yield the graded
commuting rules
\begin{eqnarray*}
\left[\partial_h^{+j},x_kT_h^{-k}\right]=\delta_{jk}I &\mbox{and}&
\left[\partial_h^{-j},x_kT_h^{+k}\right]=\delta_{jk}I,
\end{eqnarray*}

The combination of these relations with the mutual commutativity of $x_jT_h^{-j}$ resp.
$x_jT_h^{+j}$ naturally lead to the set of Weyl-Heisenberg relations (\ref{WeylHeisenbergLM}).

Based on this correspondence, the operators $X_{h}$ and $X_{-h}$ defined via the coordinate formula
\begin{eqnarray}
\label{Xhpm}X_{\varepsilon}: \f(x)\mapsto \sum_{j=1} \e_j x_j
\f(x-\varepsilon \e_j)
\end{eqnarray}
are thus the corresponding {\it Fourier duals}\footnote{In the language of Clifford
analysis it is also common to use the terminology {\it Fischer
duality} (see e.g.~\cite{Lavicka12}) instead of the {\it Fourier
duality} terminology arising in the setting of invariant theory
(cf. \cite{Howe89TransAMS}).} of $D_h^+$ and $D_h^-$,
respectively.

\subsection{Quasi-monomiality formulation on the lattice}\label{QuasiMonomialitySub}

Without loss of generality one will
consider the forward differences $\partial_h^{+j}$ as the
corresponding lowering operators embody in (\ref{quasiMonomialOp})
and the finite difference Dirac operator $D_h^+$ of forward type as the corresponding
discretization of $D=\sum_{j=1}^n \e_j \partial_{x_j}$.
From the
definitions of $\partial_h^{+j}$ and $D_h^{+}$ labeled by
(\ref{DiffPmj}) and (\ref{Dhpm}), respectively, it follows that
the multi-variable {\bf EGF}
\begin{eqnarray}
\label{factorialEGF}G_h(x,y)=\prod_{j=1}^n \left( 1+h
y_j\right)^{\frac{x_j}{h}}
\end{eqnarray}
converges asymptotically to $G(x,y)=\exp(x \cdot y)$, as $h
$ approaches to zero, and satisfies for each $y=\sum_{j=1}^ny_j\e_j\in \cl_{0,n}$
the eigenvalue equation $D_h^+ G_h(x,y)=y G_h(x,y).$

Now let us denote by $|\alpha|=\sum_{j=1}^n \alpha_j$ the sum of all the components of the  
multi-index representation $\alpha=\sum_{j=1}^n\alpha_j \e_j$ and by $\alpha!=\alpha_1!\alpha_2! \ldots \alpha_n!$ 
the corresponding multi-index factorial. 
In terms of the multi-index falling factorials of degree $|\alpha|$, defined as
\begin{eqnarray}
\label{factorialPolynomials} (x;h)_\alpha= \prod_{j=1}^n
\prod_{k=0}^{\alpha_j-1} (x_j - k h),
\end{eqnarray}
the multi-variable {\bf EGF}~(\ref{factorialEGF}) may be written as
\begin{eqnarray*}
\label{factorialEGFhypergeometric} G_h(x,y)=\prod_{j=1}^n
\exp\left( \frac{x_j}{h}\log(1+h y_j)\right) =\sum_{k=0}^\infty
\sum_{|\alpha|=k} (x;h)_{\alpha} \frac{y^\alpha}{\alpha!}.
\end{eqnarray*}

Based on the operational identity
\begin{eqnarray*}
(x;h)_\alpha=\prod_{j=1}^n \left( x_j T_h^{-j}\right)^{\alpha_j}1
\end{eqnarray*}
one can also conclude that the multi-index falling factorials
(\ref{factorialPolynomials}) are the simplest quasi-monomial
counterparts of the classical multi-variable monomials
$x^\alpha=x_1^{\alpha_1}x_2^{\alpha_2}\ldots x_n^{\alpha_n}$ on
the lattice $h\BZ^n$.
On the other hand, the Taylor series expansion
$$
\f(x+h \e_j)=\sum_{k=0}^\infty \frac{1}{k!}\left(\partial_{x_j}\right)^k\f(x)
$$
together with \cite[Theorem 2.1.1]{Roman84} gives rise
to the operational representation $T_h^{+j}=\exp\left(h\partial_{x_j}\right)$ at the level
of $\mathcal{P}=\BR[x]\otimes \cl_{0,n}$, and moreover, 
to the formal inversion formula
\begin{eqnarray}
\label{TaylorFormalInversion}
\partial_{x_j}= \frac{1}{h}\log\left( 1+h\partial_h^{+j}\right)
\end{eqnarray}
written in terms of the logarithmic function $\log(t)=\int_{1}^t
\frac{ds}{s}$. 
This in turn shows that the raising/lowering
operators considered above are a particular example of
(\ref{LMSheffer}). The functions $\kappa(t)$ and $\ell(t)$ are then
given by $\kappa(t)=1$ and $\ell(t)=\frac{\exp(ht)-1}{h}.$

\begin{remark}
When the forward finite difference operators $\partial_{h}^{+j}$ are
replaced by the backward finite difference operators $\partial_{h}^{-j}$,
it can be easily seen from (\ref{factorialEGF}) resp.
(\ref{factorialPolynomials}) that $G_{-h}(x,y)$ resp.
$\{(x;-h)_{\alpha}~:~\alpha \in \BN_0^n\}$ is the corresponding
multi-variable {\bf EGF} resp. Sheffer set of multi-variable polynomials that encodes, for each $j=1,2,\ldots,n$,
the set of raising operators $x_jT_h^{+j}$.
\end{remark}

Generally speaking, the quasi-monomials on the lattice $h\BZ^n$ constructed from the
quasi-monomial operational representation
\begin{eqnarray}
\label{quasiMonomialQuantumFieldLemma}\m_\alpha(x)=\prod_{j=1}^n
\left(M_j\right)^{\alpha_j}1
\end{eqnarray}
enables us to compute most of the families of polynomials already considered in 
\cite{SmirnovTurbiner95},~\cite{BLR98},~\cite{Kisil02},
\cite{BenCheikhZagh03},~\cite{LTW04},~\cite{BlasDattoliHorzPenson06} and \cite{VinetZhedanov09}.
In particular, based on the formal series representation 
$\partial_{x_j}=\frac{1}{h}\log T_h^{+j}$ obtained from (\ref{TaylorFormalInversion})
and on the power series expansion for $\kappa(t)$, described as
\begin{eqnarray*}
\kappa(t)=\sum_{k=0}^\infty a_k \dfrac{t^k}{k!} &\mbox{with}&\kappa(0)=a_0\neq0,
\end{eqnarray*}
one can see that the multi-variable {\bf
EGF}~$G_h(x,y;\kappa)$ defined {\it viz} equation (\ref{generalEGFhypergeometric}) encodes the set of
ladder operators 
\begin{center}
$L_j=\partial_h^{+j}$ and 
$M_j=\left(x_j-\kappa'(\partial_{x_j})\kappa\left(\partial_{x_j}\right)^{-1} \right)T_h^{- j}$.
\end{center}
Moreover, they correspond to the generators of the Weyl-Heisenberg algebra with dimension $2n+1$.

We finish this subsection with three
examples that illustrates the applicability of the {\it
quasi-monomiality} approach on the lattice $h\BZ^n$. The first two examples
involve the multi-variable Poisson-Charlier polynomials and
the multi-variable Bernoulli polynomials of the second kind. On the third example, one will sketch how the 
multi-variable quasi-monomials encoded by the central finite
difference operators $L_j=\frac{1}{2}\left(
\partial_h^{+j}+\partial_{h}^{-j}\right)$ may be computed from the multi-variable {\bf EGF} (\ref{factorialEGF}).

\begin{example}[Multi-variable Poisson-Charlier polynomials]\label{PoissonCharlierExample}
For a parameter $a \in \BR$, the multi-variable Poisson-Charlier polynomials
${\bf c}_\alpha(x;h,a)$ may be constructed from the multi-variable {\bf EGF}
$$G_h(x,y;\kappa)=\prod_{j=1}^n \exp(-ay_j)\left(
1+hy_j\right)^{\frac{x_j}{h}}.$$ 

The power series expansion $\kappa(t)$ satisfying $\kappa\left(\frac{1}{h}\log(1+ht)\right)=\exp(at)$ is given by
the function
$$\kappa(t)=\exp\left(\frac{a}{h}\left(\exp(ht)-1\right) \right).$$

Therefore, the raising operators $M_j=x_jT_h^{-j}-aI$ are thus obtained from the 
substitutions $t\rightarrow \partial_{x_j}$ on
$\frac{\kappa'(t)}{\kappa(t)}=a\exp(ht)$.
\end{example}

\begin{example}[Multi-variable Bernoulli polynomials of the second kind]\label{BernoulliExample}
Based on the function $\kappa(t)=\frac{ht}{\exp(ht)-1}$ and on the series expansion
$$\frac{1}{t}=\sum_{m=0}^\infty (1-t)^m
=\sum_{m=0}^\infty \sum_{k=0}^m \left(\begin{array}{cc}
m \\
k
\end{array}\right)(-1)^k t^k,
$$
with $0<t<2$, one obtains 
the multi-variable {\bf EGF}, logarithmic derivative and raising operators underlying to 
the multi-variable Bernoulli polynomials $\b_\alpha(x;h)$ of the second
kind. They are described as follows:
\begin{itemize}
  \item {\bf Multi-variable EGF:}
~$$G_h(x,y;\kappa)
=\prod_{j=1}^n \frac{\log\left(1+hy_j\right)}{hy_j}
~\left(1+hy_j\right)^{\frac{x_j}{h}}.$$ 
\item {\bf Logarithmic
derivative:}
$\frac{\kappa'(t)}{\kappa(t)}=\frac{1}{t}-\frac{h}{1-\exp(-ht)}$.
 \item {\bf Raising operators:}
 $$M_j=x_jT_h^{-j}-
\sum_{m=0}^\infty \sum_{k=0}^m \left(\begin{array}{ll}
m \\
k
\end{array}\right)(-1)^k \left(\frac{1}{h}\log
T_h^{+j}\right)^k-h\left(T_h^{-j}\right)^m.
 $$
\end{itemize}

Herewith $\left(\frac{1}{h}\log T_h^{+j}\right)^k$ is nothing else
than the formal representation of the iterated partial derivative $\left(\partial_{x_j}\right)^k$
on the lattice $h\BZ^n$ whereas
$\left(T_h^{-j}\right)^m=\exp\left(-mh\partial_{x_j}\right)$.
\end{example}

\begin{example}[Quasi-monomials carrying central
differences]\label{centralDifferencesExample}

When the forward finite difference operators $\partial_h^{+j}$ are replaced by the
central finite difference operators $L_j=\frac{1}{2}\left(
\partial_h^{+j}+\partial_h^{-j}\right)
$ 
one gets from its formal power series representation
$L_j=\dfrac{
1}{h}\sinh\left(h\partial_{x_j}\right)$
and from \cite[Theorem 2.1.1]{Roman84} that the multi-variable {\bf EGF} of the form
$$\prod_{j=1}^n
\exp\left(
\frac{x_j}{h}\sinh^{-1}(hy_j)\right)=\prod_{j=1}^n\left( h
y_j+\sqrt{1+h^2y_j^2}\right)^{\frac{x_j}{h}}.$$
encodes the set of ladder operators $L_j$ and $M_j=x_j[L_j,x_j]^{-1}$, with $[L_j,x_j]=\cosh(h\partial_{x_j})$.
Here one recall that the right hand side of the above formula follows from the fact
that $\sinh^{-1}(t)=\log\left( t+\sqrt{1+t^2}\right)$.

In terms of the vector-field $u: y \mapsto u(y)$, defined componentwise via 
the set of transformations
$$u_j(y)=y_j-\frac{1}{h}+\sqrt{\frac{1}{h^2}+y_j^2}$$
the above multi-variable {\bf EGF} is thus equal to $G_h(x,u(y))$,
where $G_h(x,y)$ stands the multi-variable {\bf EGF}
(\ref{factorialEGF}).

This enables to compute the quasi-monomials $\m_\alpha(x)$ generated from the
operational rule (\ref{quasiMonomialQuantumFieldLemma}) through the action of the multi-index derivative
$$\left(\partial_y\right)^\alpha=
\left(\partial_{y_1}\right)^{\alpha_1}\left(\partial_{y_2}\right)^{\alpha_2}\ldots
\left(\partial_{y_n}\right)^{\alpha_n}$$ on $G_h(x,u(y))$. In
concrete, from (\ref{factorialEGFhypergeometric}) one gets
\begin{eqnarray*}
\m_\alpha(x)=\left[\left(\partial_y\right)^\alpha
G_h(x,u(y))\right]_{y={\bf
0}}=\frac{\gamma_\alpha}{\alpha!}~(x;h)_\alpha
\end{eqnarray*}
with $\gamma_\alpha=\left[ (\partial_y)^\alpha
u(y)^\alpha\right]_{y={\bf 0}}$. 

Here one recall that $(x;h)_\alpha$ denotes the multi-index falling
factorial of degree $|\alpha|$ defined {\it viz} (\ref{factorialEGF}).
\end{example}

\begin{remark}
In contrast with \cite[Example 4 of Section 2]{DHS96}, where the
quasi-monomials carrying central finite difference operators were
computed by a binomial convolution formula, the quasi-monomials obtained
in Example \ref{centralDifferencesExample}
correspond to the 'Taylor
coefficients' of the multi-variable {\bf EGF} $G_h(x,u(y))$.
\end{remark}

\begin{remark}
The Bernoulli polynomials of the second kind considered in Example \ref{BernoulliExample} may also be 
formulated via formal series representations of integral operators.

One suggest the interested reader to take a look for the sequence of examples 
explored in \cite[Subsection 4.2]{BLR98}, on which families of Bernoulli polynomials 
of the second kind, encoded by a central finite difference operator, were described 
in terms of the Bessel functions $
J_s(u)=\frac{1}{\Gamma(s+1)}\left(\frac{u}{2}\right)^{s}{~}_0 F_1\left(s+1;-\frac{u^2}{4}\right)
$.
\end{remark}

\section{The hypercomplex approach}\label{HypercomplexApproach}

\subsection{Classes of Clifford-vector-valued raising operators}

Based on the Lie-algebraic description for the 
Fischer duals of $D_h^+$, we now proceed to the construction of hypercomplex extensions for the quasi-monomial basis 
(\ref{quasiMonomialQuantumFieldLemma}). Such description, as obtained in the following proposition,
encompasses the set of ladder operators $M_j$ obtained in Subsection \ref{QuasiMonomialitySub}.

\begin{proposition}\label{FourierDualProp}
Let $\kappa(t)$ defined as above and $X_h$ the multiplication operator defined
{\it viz} (\ref{Xhpm}). If there is a multi-variable function
$\lambda(y)$ ($y\in \BR^n$) such that
$$\lambda\left( \frac{D_h^+\exp(x\cdot y)}{\exp(x\cdot y)}\right)=\prod_{j=1}^n\kappa(y_j),$$
then the Fourier dual $\Lambda_h$ of $D_h^+$ is given by
\begin{eqnarray*}
 \Lambda_h&=&
X_h-\left[\log \lambda\left(D_h^+\right),x\right].
\end{eqnarray*}

\end{proposition}

\proof First, recall that from \cite[Theorem 3.6.5]{Roman84} and
from the isomorphism between the algebra of formal power series
and the algebra of linear functionals (cf. \cite[Theorem
2.1.1]{Roman84}) there is
a one-to-one correspondence between the logarithmic derivative
$$\frac{\kappa'(y_j)}{\kappa(y_j)}=\left[\frac{d
\log(\kappa(t))}{dt}\right]_{t=y_j}$$ and the Pincherle derivative
$\left[\log\left(\kappa\left(\partial_{x_j}\right)\right),x_j\right]
=\kappa'(\partial_{x_j})\kappa(\partial_{x_j})^{-1}$.

In the same order of ideas, for $L_j=\partial_h^{+j}$ there is a
one-to-one correspondence between the backward shift
$T_h^{-j}=\left(T_h^{+j}\right)^{-1}=\left[\partial_h^{+j},x_j\right]^{-1}$
that yields from the product rules (\ref{productRule}) and the
exponentiation relation $\exp(-hy_j)=\exp(hy_j)^{-1}$.

Then, the {\it Fourier dual} $\Lambda_h=\sum_{j=1}^n \e_jM_j$
constructed from (\ref{LMSheffer}) have the following
Lie-algebraic representation in the algebra
$\mbox{End}(\mathcal{P})$:
\begin{eqnarray}
\label{FourierDualDh+} \Lambda_h=X_h-\sum_{j=1}^n \e_j
\kappa'(\partial_{x_j})\kappa(\partial_{x_j})^{-1}T_h^{-j}.
\end{eqnarray}

 Now let us take a close look for the commutator
$\left[\log\lambda(D_h^+),x\right]$. A short computation based on
the identity $\exp\left(\left(x+h\e_j\right)\cdot
y\right)=\exp(x\cdot y)\exp(hy_j)$ shows that the quantity
$$\frac{D_h^+\exp(x\cdot y)}{\exp(x\cdot y)}=\sum_{j=1}^n \e_j \frac{\exp(hy_j)-1}{h}$$
corresponds to the representation of $D_h^+$ in the algebra of formal power
series expansions (cf.~\cite[Theorem
2.1.1]{Roman84}).

Combination of the chain rule
$$\partial_{y_j}(\log \lambda)\left(\frac{D_h^+\exp(x\cdot
y)}{\exp(x\cdot y)}
\right)~\exp(hy_j)=\frac{\kappa'(y_j)}{\kappa(y_j)}$$ with
\cite[Theorem 3.6.5]{Roman84} result into the sequence of
identities
\begin{eqnarray*}
\left[\log \lambda(D_h^+),x\right]&=&\sum_{j=1}^n \e_j
\left[\log\lambda(D_h^+),x_j\right] \\
&=&\sum_{j=1}^n \e_j \left[
\log\left(\kappa\left(\partial_{x_j}\right)\right),x_j\right]
\exp(-h\partial_{x_j}) \\
&=& \sum_{j=1}^n \e_j
\kappa'(\partial_{x_j})\kappa(\partial_{x_j})^{-1}T_h^{-j}.
\end{eqnarray*}

Hence, the equation (\ref{FourierDualDh+}) is equivalent to
$$\Lambda_h=X_h-\left[\log\lambda\left( D_h^+\right),x\right].$$
\qed

Recall that from Proposition \ref{FourierDualProp}, the set of formal inversion formulae 
(\ref{TaylorFormalInversion})
together with \cite[Theorem 2.1.1]{Roman84} even shows that the multi-variable function $\lambda(y)$ ($y \in \BR^n$)
always exists and it is explicitly given by
\begin{eqnarray}
\label{LambdaY}\lambda(y)=\prod_{j=1}^n \kappa\left(\frac{1}{h}\log(1+hy_j)\right).
\end{eqnarray}

So, one looks for the ladder operators $\Lambda_h$ as covariant versions of the
$X_h$ defined {\it viz} equation (\ref{Xhpm}) on which the operator $\lambda(D_h^+)$
is obtained from the substitutions $y_j\rightarrow \partial_h^{+j}$ on the right hand side 
of (\ref{LambdaY}). 

Regardless this construction one can easily see that
the multi-variable function of the form
$\lambda(y)=\prod_{j=1}^n  (1+hy_j)^{d_j}$
($d_j\in \BR$) yields $\Lambda_h=\sum_{j=1}^n
\e_j(x_j-d_jh)T_h^{-j}$ as the finite difference
counterpart of the multiplication operator
$X=\sum_{j=1}^n\e_j~x_jI$ on $h\BZ^n$.

Analogously, one can also use the same scheme to construct finite
difference discretizations for the Clifford-Hermite operator
(cf.~\cite[Subsection 3.2]{CFK11})
\begin{eqnarray}
\label{CliffordHermiteOp}
X-D=-\exp\left(\frac{|x|^2}{2}\right)D\exp\left(-\frac{|x|^2}{2}\right).
\end{eqnarray}

For example, the Fischer dual $\Lambda_h=X_h-D_h^-$ of $D_h^+$, with $D_h^-=\left[\log\lambda(D_h^+),x\right]$, 
arises in case where $\lambda(y)$ satisfies the $\log-$constraint $$\log \lambda\left(
\frac{D_h^+\exp(x\cdot y)}{\exp(x\cdot y)}\right)=\sum_{j=1}^n
 \frac{\exp(hy_j)-hy_j}{h^2}.$$ 

As illustrated on the next example, such operational scheme
enable us a fully rigorous way to obtain a wide class of Fischer duals for $D_h^+$.

\begin{example}[A non-trivial discretization for the Clifford-Hermite operator]
For the multi-variable function $\lambda(y)$ defined as  
$$\lambda(y)=\prod_{j=1}^n \exp\left(\dfrac{1+hy_j}{2h^2}+\dfrac{1}{(2+2hy_j)h^2}\right)$$
one has 
$$
\lambda\left(\frac{D_h^+ \exp\left(x\cdot y\right)}{\exp\left(x\cdot y\right)}\right)
=\prod_{j=1}^n \exp\left(\frac{\cosh(hy_j)}{h^2}\right).
$$

Moreover, the set of identities $(\partial_{y_j}\cosh)(hy_j)=h\sinh(hy_j)$ and
$$\frac{\sinh(hy_j)}{h}\exp(-hy_j)=\frac{1-\exp(-2hy_j)}{2h}$$
leads to $\left[\log \lambda(D_h^+),x_j\right]=\partial_{2h}^{-j}$, and hence, to
$$
[\log \lambda(D_h^+),x]=\sum_{j=1}^n\e_j [\log \lambda(D_h^+),x_j]=D_{2h}^-.
$$

Under the above choice for $\lambda(y)$, the discretization for the Clifford-Hermite operator
(\ref{CliffordHermiteOp})
is thus given by $\Lambda_h=X_h-D_{2h}^-$, where $D_{2h}^-$
stands the finite difference Dirac operator of backward type
underlying to the coarse lattice $(2h)\BZ^n$ of $h\BZ^n$.
\end{example}

\subsection{Appell set formulation for finite difference Dirac operators}\label{AppellSetFormulation}

Once developed over the previous subsection the key tools to formulate the hypercomplex extension
of the quasi-monomiality principle (\ref{quasiMonomialOp}), one have now the minimal amount of tools required
to construct Appell sets of Clifford-vector-valued polynomials on the lattice $h\BZ^n$.

One say that the set $\{ \w_k(x;h;\lambda)~:~k\in
\BN_0\}$ of
Clifford-vector-valued polynomials is an Appell set with respect to $D_h^+$ if $
\w_0(x;h;\lambda)=\a$ is a Clifford number and $D_h^+
\w_k(x;h;\lambda)$ is a Clifford-vector-valued polynomial of
degree $k-1$ satisfying the Appell set property
\begin{eqnarray}
\label{AppellDh+}D_h^+ \w_k(x;h;\lambda)=k\w_{k-1}(x;h;\lambda).
\end{eqnarray}

Equivalently, the construction of Appell sets may be
formulated as a time-evolution problem in the space-time
domain $ h\BZ^n \times \BR$, described as follows: Find for each
$(x,t)\in h\BZ^n \times \BR$ a function
$\G_h(x,t;\lambda)$ that satisfies the set of equations
\begin{eqnarray}
\label{evolutionEGF} \left\{\begin{array}{lll} D_h^+
\G_h(x,t;\lambda)=t\G_h(x,t;\lambda) & \mbox{for} & (x,t)\in
h\BZ^n \times \BR \setminus \{0\}
\\ \ \\
\G_h(x,0;\lambda)=\a & \mbox{for} & x\in h\BZ^n 
.
\end{array}\right.
\end{eqnarray}

Based on the embedding of $(x,t) \in \BR^{n+1}$ in
$\cl_{0,n}$ through the paravector representation
$t+x=t+\sum_{j=1}^n \e_j x_j$, one may determine, in the same order
of ideas of \cite[Section 5]{CFM11}, the solution of
(\ref{evolutionEGF}) as an hypercomplex version of the Taylor
series expansion.
Indeed, if $\G_h(x,t;\lambda)$ is a $C^\infty$-function with
respect to $t\in \BR$ such that the $k-$derivative term
$$\left[ \left(\partial_t\right)^k\G_h(x,t;\lambda)\right]_{t=0}=\w_k(x;h;\lambda)$$
is a Clifford-vector-valued polynomial of degree $k$, then the Taylor
series expansion of $\G_h(x,t;\lambda)$ around $x \in h\BZ^n$ given by
\begin{eqnarray}
\label{EGFDh+}\G_h(x,t;\lambda)=\sum_{k=0}^\infty
\frac{t^k}{k!}~\w_k(x;h;\lambda)
\end{eqnarray}
{\it uniquely} determines the hypercomplex {\bf EGF} for $D_h^+$ as solution
of (\ref{evolutionEGF}). 

Now let us recast the quasi-monomiality principle in terms of the Fischer dual $\Lambda_h$.
Based on the Fock space formalism, each Clifford-vector-valued polynomial 
$\w_k(x;h;\lambda)$ is constructed by means of the operational
rule
\begin{eqnarray}
\label{AppellDiscrete} \w_k(x;h;\lambda)=\mu_k\left(
\Lambda_h\right)^k\a.
\end{eqnarray}

 The constants $\mu_k\in \BR$ ($k \in \BN_0$) given in (\ref{AppellDiscrete}) 
 are thus determined from the condition $\w_0(x;h;\lambda)=\a$ (normalization condition) and from the Appell set
 constraint (\ref{AppellDh+}). Such operational representation enables to compute the hypercomplex 
 {\bf EGF} (\ref{EGFDh+}) 
 through the operational representation $\G_h(x,t;\lambda)=\G(\Lambda_h,t)\a$, where the function
\begin{eqnarray}
\label{EGFDirac} \G(x,t)=\sum_{k=0}^\infty \frac{t^k}{k!}~\mu_kx^k
\end{eqnarray}
corresponds to the hypercomplex {\bf EGF} encoded by the {\it continuum} Dirac operator
$D=\sum_{j=1}^n \e_j \partial_{x_j}$.

The bold notations $\G$ and $\G_h$ are adopted throughout to make a clear distinction between the
hypercomplex {\bf EGF}
of the above form and the multi-variable {\bf EGF} $G$ and $G_h$, respectively, already introduced in Subsection
\ref{QuasiMonomialitySub}. In addition, the shortland notation $\w_k(x;h)$ resp. $\G_h(x,t)$ will be used
when one refers to the Appell polynomial of degree $k$ resp. to the hypercomplex {\bf EGF}
determined from the constant function $\lambda(y)=1$
, that is, the Clifford-vector-valued polynomials resp. the hypercomplex {\bf EGF}
generated from the Fischer dual $X_h$ defined in (\ref{Xhpm}).

The next example, corresponding to hypercomplex extension of the multi-index
falling factorials (\ref{factorialPolynomials}), will be of special interest in the forthcoming subsections.

\begin{example}[The hypercomplex version of the multi-index falling factorials]
\label{HypercomplexFactorialPolynomials}
 The hypercomplex extension of the multi-index falling factorials $(x;h)_\alpha$ of order $|\alpha|=k$ provided by
 (\ref{factorialPolynomials}) is represented, for a given $\a \in \cl_{0,n}$, 
 through the operational formula $$\w_k(x;h)=\mu_k \left(X_h\right)^k \a.$$

From the Weyl-Heisenberg graded commuting relations (\ref{WeylHeisenbergLM}) encoded by the set of ladder operators 
$M_j=x_jT_h^{-j}$, is it clear that $\left(X_h\right)^2$ is a scalar-valued operator given by the summation
formula
$$
\left(X_h\right)^2=-\sum_{j=1}^n \left( x_jT_h^{-j}\right)^2.
$$ 

A straightforward computation
based on the multinomial formula gives rise to
\begin{eqnarray*}
\w_{2m}(x;h) &=&\mu_{2m}\left( -\sum_{j=1}^n
  \left( x_jT_h^{-j}\right)^2\right)^{m}\a \\
&=&(-1)^m\mu_{2m}\sum_{s=0}^m\sum_{|\alpha|=s} \frac{m!}{\alpha!}
~(x;h)_{2\alpha}~\a
\\ \ \\
\w_{2m+1}(x;h)&=&\frac{\mu_{2m+1}}{\mu_{2m}}\sum_{j=1}^n
\e_j~x_jT_h^{-j}\w_{2m}(x;h)
\\
&=& (-1)^m\mu_{2m+1}\sum_{j=1}^n\sum_{s=0}^m\sum_{|\alpha|=s}
\frac{m!}{\alpha!}~(x;h)_{2\alpha+\e_j}~\e_j\a.
\end{eqnarray*}

Next, one will recast the above Clifford-vector-valued polynomials
in terms of the inverse of the $n-$dimensional Weierstra\ss~transform
$$
(W \f)(x)=\frac{1}{(4\pi)^{\frac{n}{2}}}\int_{\BR^n} \f(x-y)\exp\left(-\frac{|y|^2}{4}\right)dy=
\exp\left(-\frac{1}{2}D^2\right)\f(x).
$$

Here one recalls that $\exp\left(-\frac{1}{2}D^2\right)
=\prod_{j=1}^n\exp\left(\frac{1}{2}\left(\partial_{x_j}\right)^2\right)$
corresponds the operational representation of the $n-$dimensional Weierstra\ss~operator $W$ 
on $\mathcal{P}=\BR[x]\otimes \cl_{0,n}$ (cf.~\cite[Section 5.3]{BLR98}).
In particular, each multi-variable Hermite polynomial
$H_\beta(x)$ of order $|\beta|$, determined from 
the action of $W^{-1}$ on the multi-index polynomial $
x^\beta=x_1^{\beta_1}x_2^{\beta_2}\ldots x_n^{\beta_n}$ of degree $|\beta|$, is thereby represented as
\begin{eqnarray*}
H_{\beta}(x)&=&\exp\left(\frac{1}{2}D^2\right)x^\beta
\\&=&
\sum_{k=0}^{\left\lfloor \frac{|\beta|}{2} \right\rfloor}  
\sum_{|\alpha|=k}\frac{1}{\alpha!}\left( -\frac{1}{2h^2}\right)^{|\alpha|}(\beta h;h)_{2\alpha}~
x^{\beta-\alpha}.
\end{eqnarray*}

Now set $\e=\sum_{j=1}^n\e_j$ and $|x|_1=\sum_{j=1}^n \e_j |x_j|$. 
From the constraint $\sum_{j=1}^n |x_j|=2mh$ 
and the identity
$$
\sum_{s=0}^m\sum_{|\alpha|=s} \frac{m!}{\alpha!}
~(x;h)_{2\alpha}= m!~\left(2h^2\right)^{2m}
H_{\frac{|x|_1}{h}}\left(-\frac{1}{2h^2}\e\right)
$$
one gets the following closed formula for $\w_{2m}(x;h)$ on the lattice $h\BZ^n$, written in terms of
$W^{-1}=\exp\left(\frac{1}{2}D^2\right)$:
$$\w_{2m}(x;h)=(-1)^m\mu_{2m}m!~\left(2h^2\right)^{2m} 
\left(W^{-1}y^{\frac{|x|_1}{h}}\right)\left(-\frac{1}{2h^2}\e\right)\a.$$

Moreover, a closed formula for $\w_{2m+1}(x;h)$ follows straightforwardly from
combination of the above expression with the recursive relation
$$\w_{2m+1}(x;h)=\frac{\mu_{2m+1}}{\mu_{2m}}\w_{2m}(x;h).$$
\end{example}

\subsection{Classes of discrete hypercomplex polynomials}

Next, one will look further for the construction of new classes of hypercomplex polynomials of discrete
variable based on a Rodrigues-type representation of the
operational formula (\ref{AppellDiscrete}). That consists into 
the construction of an operator $\sigma(D_h^+)$ that intertwines the Fischer duals 
$\Lambda_h=X_h-\left[ \log \lambda(D_h^+),x\right]$ and $X_h$ of $D_h^+$.

The subsequent lemma provides us the required ingredient to derive such representation formula in
Proposition \ref{RodriguesFormulaXh}.

\begin{lemma}\label{AexpB}
If the three generators $A,B,C$ of a certain Lie algebra satisfy the
graded commutation rules $[A,B]=C$ and $[C,B]=0$, then for the
exponentiation operator
$$ \exp(B)=\sum_{k=0}^\infty \frac{1}{k!}B^k$$
one gets $[A,\exp(B)]=C  \exp(B)$.
\end{lemma}

\proof First, recall the following summation formula that holds for
every $A,B$ and $k \in \BN$:
\begin{eqnarray*}
\left[A,B^k\right]=\sum_{j=0}^{k-1}B^{j}[A,B] B^{k-1-j}.
\end{eqnarray*}

Under the conditions $[A,B]=C$ and $[C,B]=0$ one gets that the above
summation formula equals to $\left[A,B^k\right]=C B^{k-1}$. This
leads to
\begin{eqnarray*}
[A,\exp(B)]&=&\sum_{k=0}^\infty\frac{1}{k!}\left[A,B^k\right] \\
&=& \sum_{k=1}^\infty\frac{1}{(k-1)!}CB^{k-1} \\
&=&C \exp(B).
\end{eqnarray*}
\qed

\begin{proposition}[Rodrigues-type formula]\label{RodriguesFormulaXh}
Let $\lambda(y)$ and $\Lambda_h=X_h-\left[ \log \lambda(D_h^+),x\right]$ be the multi-variable function and
the Fischer dual of $D_h^+$ given by Proposition
\ref{FourierDualProp}, respectively, and $\{\w_k(x;h)~:~k\in\BN_0\}$ the Appell set of
Clifford-vector-valued polynomials obtained in Example \ref{HypercomplexFactorialPolynomials}.

Then, the operator $\sigma\left(D_h^+\right)\in
\mbox{End}(\mathcal{P})$ defined {\it viz}
$$\sigma(D_h^+)=\lambda\left[
\exp(-x\cdot D)D_h^+\exp(x\cdot D)\right],$$ with $x\cdot D=\sum_{j=1}^n x_j \partial_{x_j}$, satisfies
$\Lambda_h=\sigma(D_h^+)^{-1}~X_h~\sigma(D_h^+).$

Moreover, the Clifford-vector-valued polynomials $\w_k(x;h;\lambda)$, 
determined from the operational formula (\ref{AppellDiscrete}), are given by
  $$ \w_k(x;h;\lambda)=\lambda({\bf 0})~\sigma(D_h^+)^{-1}\left[\w_k(x;h)\right].$$
\end{proposition}

\proof Under the conditions of Proposition \ref{FourierDualProp}
one can see that the multi-variable function $\lambda(y)$
satisfies
$$ \lambda\left( \frac{D_h^+ \exp(x\cdot y)}{\exp(x\cdot y)}\right)
=\prod_{j=1}^n  \kappa(y_j).$$

Thus, under the substitutions $y_j \rightarrow \partial_{x_j}$ on both sides
of the above equation, one immediately gets that $\sigma(D_h^+)=\lambda\left[ \exp(-x\cdot
D)D_h^+\exp(x\cdot D)\right]$ is shift-invariant and satisfies the log-relation
$$ \log \sigma (D_h^+)=\sum_{j=1}^n \log \kappa\left(\partial_{x_j}\right).$$

Based on the same order of ideas considered in the proof of
Proposition \ref{FourierDualProp}, the following set of graded commuting relations
 $$[x_j,-\log \sigma (D_h^+)]=\left[\log \kappa\left(\partial_{x_j}\right),x_j \right]
 =\kappa'\left(\partial_{x_j}\right)\kappa\left(\partial_{x_j}\right)^{-1}$$
  follows straightforwardly for each $j=1,2,\ldots,n$.
  
A direct application of Lemma \ref{AexpB} through the substitutions
$A=x_jI$, $B=\exp\left(-\log \sigma (D_h^+)\right)$ and
$C=\kappa'\left(\partial_{x_j}\right)\kappa\left(\partial_{x_j}\right)^{-1}$
gives rise to
$$\left[x_j,\sigma(D_h^+)^{-1}\right]
=\left[x_j,\exp\left(-\log \sigma(D_h^+)\right)\right]
=\kappa'\left(\partial_{x_j}\right)\kappa\left(\partial_{x_j}\right)^{-1}
\sigma(D_h^+)^{-1}.
$$

By multiplying both sides of the above identity on the right by $T_h^{-j}=\exp(-h\partial_{x_j})$
one obtains, after
a straightforward computation, the following set of intertwining
formulae
$$
\left(x_j-\kappa'\left(\partial_{x_j}\right)\kappa\left(\partial_{x_j}\right)^{-1}\right)
T_h^{-j} \sigma (D_h^+)^{-1} = \sigma (D_h^+)^{-1}x_jT_h^{-j}.
$$

Hence, the relation
$$\left(X_h-\left[
\log\lambda(D_h^+),x\right]\right)\sigma(D_h^+)^{-1}=\sigma(D_h^+)^{-1}~X_h$$
follows straightforwardly from linearity arguments.

Finally, by multiplying both sides of the above relation on the right by $\sigma(D_h^+)$ results
into the operational representation $\Lambda_h=\sigma(D_h^+)^{-1}~X_h~\sigma(D_h^+).$

The statement for $\w_k(x;h;\lambda)=\lambda({\bf 0})\sigma(D_h^+)^{-1}\left[\w_k(x;h)\right]$ is then
immediate from the previous operational representation formula for $\Lambda_h$ and from the basic relations
$\lambda({\bf 0})=\kappa(0)^n$ (${\bf 0}\in\BR^n$) and
$$\prod_{j=1}^n\kappa\left(\partial_{x_j}\right)\a=\prod_{j=1}^n \kappa(0)\a=\kappa(0)^n\a.$$ 
\qed

Based on the integral representation 
$$
\sigma\left(D_h^+\right)^{-1}=\int_{0}^\infty \exp\left(-\sigma(D_h^+)s\right)~ds
$$ 
resulting from the substitution $r \rightarrow \sigma\left(D_h^+\right)$ on the identity
$r^{-1}=\int_{0}^\infty \exp(-rs)ds$,
the next corollary is rather obvious.

\begin{corollary}[Integral representation]\label{RodriguesFormulaXhCorollary}
Under the conditions of Theorem \ref{RodriguesFormulaXh}, each Clifford-vector-valued polynomial $\w_k(x;h;\lambda)$
admits the following integral representation formula
$$
\w_k(x;h;\lambda)=\lambda({\bf 0})\int_{0}^\infty \exp\left(-s\sigma(D_h^+)\right)\left[\w_k(x;h)\right]ds.
$$
\end{corollary}

\begin{example}[Hypercomplex extension of Poisson-Charlier polynomials]\label{HypercomplexPoissonCharlier}
Let us consider again the function
$\kappa(t)=\exp\left(\frac{a}{h}\left(\exp(ht)-1\right)\right)$ used
in Example \ref{PoissonCharlierExample} to characterize the
multi-variable Poisson-Charlier polynomials ${\bf
c}_{\alpha}(x;h;a)$ of order $|\alpha|$.

Following the same order of ideas of Example \ref{HypercomplexFactorialPolynomials},
the Clifford-vector-valued operator of the form $\Lambda_h=\sum_{j=1}^n \e_j(x_jT_h^{-j}-aI)$
gives rise to the formulae
\begin{eqnarray*}
\w_{2m}(x;h;\lambda) &=&\mu_{2m}\left( -\sum_{j=1}^n
(x_jT_h^{-j}-a)^2\right)^{m}\a \\
&=&(-1)^m\mu_{2m}\sum_{s=0}^m\sum_{|\alpha|=s} \frac{m!}{\alpha!}
~{\bf c}_{2\alpha}(x;h;a)~\a
\\ \ \\
\w_{2m+1}(x;h;\lambda)&=&\frac{\mu_{2m+1}}{\mu_{2m}}\sum_{j=1}^n
\e_j(x_jT_h^{-j}-a)\w_{2m}(x;h;\lambda)
\\
&=& (-1)^m\mu_{2m+1}\sum_{j=1}^n\sum_{s=0}^m\sum_{|\alpha|=s}
\frac{m!}{\alpha!}~{\bf c}_{2\alpha+\e_j}(x;h;a)~\e_j\a.
\end{eqnarray*}

Alternatively, based on the properties
$\exp\left( a\partial_h^{+j}\right)\a=\a$,
$\exp\left( -a\partial_h^{+j}\right)=\exp\left( a\partial_h^{+j}\right)^{-1}$ and
$$\prod_{j=1}^n\exp\left( -a\partial_h^{+j}\right)
=\exp\left( -a\sum_{j=1}^n \partial_{h}^{+j}\right)
$$
it follows
from direct application of Proposition \ref{FourierDualProp} and
Proposition~\ref{RodriguesFormulaXh} that each $\w_{k}(x;h;\lambda)$ admits the
following operational representation
$$ \w_{k}(x;h;\lambda)=\exp\left( -a\sum_{j=1}^n \partial_{h}^{+j}\right)\w_{k}(x;h).$$
\end{example}

\begin{remark}[The Poisson-Charlier connection]\label{HypercomplexPoissonCharlierRemark}
 It is straighforward to see that one can recast $\Lambda_h=X_h-D_h^-$ as
$$\Lambda_h=\sum_{j=1}^n \e_j \left(
\left(x_j+\frac{1}{h}\right)T_h^{-j}-\frac{1}{h}I\right).$$

Thus, the polynomials $\w_k\left(x,h;\lambda\right)$ determined from the
operatorial rule (\ref{AppellDiscrete}) may be described in terms
of the hypercomplex Poisson-Charlier polynomials provided by Example
\ref{HypercomplexPoissonCharlier}. Through the parameter substitution $a=\frac{1}{h}$,
such description is obtained in terms displaced vector variable $x+\frac{1}{h}\e$, with  
$\e=\sum_{j=1}^n \e_j$.
\end{remark}

\subsection{Hypergeometric series representations}\label{HypergeometricSubsection}

One have now the key ingredients to compute a closed formula for the hypercomplex {\bf EGF} formulated in
Subsection \ref{AppellSetFormulation} from the set of equations (\ref{evolutionEGF}). 
Based on the operational formula $\G_h(x,t;\lambda)=\G(\Lambda_h,t)\a$
one shall compute, first of all, the constants $\mu_k$ assigned by the operational formula 
(\ref{AppellDiscrete}) from the constraint $$D \G(x,t)=t\G(x,t).$$ 

The next lemma, a particular case of \cite[Lemma 3.1]{FR11},
is the key ingredient used in the proof of Proposition \ref{EGFDiracProposition}. 
Recall that the Clifford-vector-valued monomials 
$x^k$ arising in the ansatz (\ref{EGFDirac}) satisfy $x^{2m}=(-1)^m|x|^{2m}$ ($k=2m$)
and $x^{2m+1}=(-1)^mx|x|^{2m}$ ($k=2m+1$).
\begin{lemma}\label{Dxk}
The Clifford-vector-valued monomials 
satisfy the recursive relations
\begin{eqnarray}
\label{MonomialityXk}   D x^{k}&=& \left\{
\begin{array}{lll}
-2m~x^{2m-1}, & \mbox{if}~~k=2m\in \BN_0 \\
-(2m+n)x^{2m} & \mbox{if}~~k=2m+1\in \BN
\end{array}\right.
\end{eqnarray}
\end{lemma}

\begin{proposition}\label{EGFDiracProposition}
The {\bf EGF} (\ref{EGFDirac}) satisfies $D \G(x,t)=t\G(x,t)$ if
and only if the constants $\mu_k$ are equal to
\begin{eqnarray*}
\mu_{2m}=(-1)^m\frac{\left(\frac{1}{2}\right)_{m}}{\left(\frac{n}{2}\right)_{m}} &\mbox{and}&
\mu_{2m+1}=(-1)^m\frac{\left(\frac{3}{2}\right)_{m}}{\left(\frac{n}{2}+1\right)_{m}}.
\end{eqnarray*}

Moreover
\begin{eqnarray*}
\G(x,t)
&=&{~}_0F_1\left( \frac{n}{2};-\frac{t^2}{4}x^2\right)
+tx{~}_0F_1\left( \frac{n}{2}+1;-\frac{t^2}{4}x^2\right).
\end{eqnarray*}
\end{proposition}

\proof 
Direct application of Lemma \ref{Dxk} to each summand of (\ref{EGFDirac}) leads to the
splitting formula
\begin{eqnarray*}
D \G(x,t)=\sum_{m=0}^\infty -2m \mu_{2m} x^{2m-1}
~\frac{t^{2m}}{(2m)!}+ \sum_{m=0}^\infty -(2m+n) \mu_{2m+1} x^{2m}
~\frac{t^{2m+1}}{(2m+1)!}.
\end{eqnarray*}

Therefore $\G(x,t)$ is a solution of the equation $D \G(x,t)=t\G(x,t)$ if and only if
\begin{eqnarray*}
\frac{-2m\mu_{2m}}{(2m)!}=\frac{\mu_{2m-1}}{(2m-1)!} &\mbox{and}&
\frac{-(2m+n)\mu_{2m+1}}{(2m+1)!}=\frac{\mu_{2m}}{(2m)!},
\end{eqnarray*}
that is
\begin{eqnarray*}
\mu_{2m}=-\mu_{2m-1} & \mbox{and} & \mu_{2m+1}=-\frac{2m+1}{2m+n}\mu_{2m-1}.
\end{eqnarray*}

Induction over $m \in \BN_0$ shows that the $\mu_k$'s are explicitly given by the formulae
\begin{center}
$\mu_{2m}=(-1)^m\frac{\left(\frac{1}{2}\right)_{m}}{\left(\frac{n}{2}\right)_{m}}$ and
$\mu_{2m+1}=(-1)^m\frac{\left(\frac{3}{2}\right)_{m}}{\left(\frac{n}{2}+1\right)_{m}}$,
\end{center}
where $(a)_m=\frac{\Gamma(a+m)}{\Gamma(a)}$ stands the Pochhammer symbol. 

Hence, a short computation based on the properties 
$$\frac{\left(\frac{1}{2}\right)_{m}}{(2m)!}
=\frac{\left(\frac{1}{4}\right)^m}{m!}=\frac{\left(\frac{3}{2}\right)_{m}}{(2m+1)!}$$ 
gives in turn the above hypergeometric series splitting of type ${~}_0F_1$ for $\G(x,t)$:
\begin{eqnarray*}
\G(x,t)&=& \sum_{m=0}^\infty
(-1)^m\frac{\left(\frac{1}{2}\right)_{m}}{\left(\frac{n}{2}\right)_{m}}
x^{2m} \frac{t^{2m}}{(2m)!}+ \sum_{m=0}^\infty
(-1)^m\frac{\left(\frac{3}{2}\right)_{m}}{\left(\frac{n}{2}+1\right)_{m}}
x^{2m+1} \frac{t^{2m+1}}{(2m+1)!}\\
&=&{~}_0F_1\left( \frac{n}{2};-\frac{t^2}{4}x^2\right)+tx{~}_0F_1\left( \frac{n}{2}+1;-\frac{t^2}{4}x^2\right).
\end{eqnarray*}
\qed

Direct combination of the above proposition with Proposition \ref{RodriguesFormulaXh}, the next
corollary follows naturally.
\begin{corollary}\label{hypercomplexEGFCorollary}
The hypercomplex {\bf EGF} $\G_h(x,t)$ defined for $\lambda(y)=1$ via the Taylor series
expansion (\ref{EGFDh+}) has the formal hypergeometric
series representation
$$
\G_h(x,t)={~}_0F_1\left( \frac{n}{2};-\frac{t^2}{4}\left(X_h\right)^2\right)\a+
tX_h{~}_0F_1\left( \frac{n}{2}+1;-\frac{t^2}{4}\left(X_h\right)^2\right)\a.
$$

Moreover, under the conditions of Proposition \ref{RodriguesFormulaXh} the hypercomplex
{\bf EGF} $\G_h(x,t;\lambda)$ is given by
$$\G_h(x,t;\lambda)=\lambda({\bf 0})\sigma(D_h^+)^{-1}\left[\G_h(x,t)\right].$$
\end{corollary}

\begin{remark}\label{BesselConnection}
The ${~}_0F_1$-hypergeometric series representation obtained
on Corollary \ref{hypercomplexEGFCorollary} shows that the hypercomplex {\bf EGF} $\G_h(x,t;\lambda)$
that arises from the substitution $x \rightarrow \Lambda_h$ on $\G(x,t)\a$ is of Bessel type and 
close to the monogenic exponential function obtained in \cite[Theorem 2]{CFM11}.
Indeed, based on the hypergeometric representation for the Bessel function $J_s(u)$ of order $s$:
$$
J_s(u)=\frac{1}{\Gamma(s+1)}\left(\frac{u}{2}\right)^{s}{~}_0 F_1\left(s+1;-\frac{u^2}{4}\right),
$$
the hypercomplex {\bf EGF} $\G(x,t)$ obtained in Proposition \ref{EGFDiracProposition} 
is equivalent to 
$$
\G(x,t)=\Gamma\left( \frac{n}{2}\right)\left(\frac{tx}{2}\right)^{-\frac{n}{2}+1}
\left(J_{\frac{n}{2}-1}(tx)+n~J_{\frac{n}{2}}(tx)\right).
$$
\end{remark}

\begin{remark}\label{hypercomplexEGFRemark}
Based on Example \ref{HypercomplexFactorialPolynomials}  
the even powers $(X_h)^{2m}\a$ are described, for each $x \in h\BZ^n$ satisfying $\|x\|_1:=\sum_{j=1}^n|x_j|=2mh$,
by the set of identities
\begin{eqnarray*}
(X_h)^{2m}\a&=&(-1)^m\frac{\left(\frac{n}{2}\right)_{m}}
{\left(\frac{1}{2}\right)_{m}}\w_{2m}(x;h)\\
&=&
(-1)^m m!\left(-\frac{1}{2h}\e\right)^{-\frac{|x|_1}{h}} 
H_{\frac{|x|_1}{h}}\left(-\frac{1}{2h}\e\right)\a,
\end{eqnarray*}
where $H_\beta(y)=\exp\left(\frac{1}{2}D^2\right)y^\beta$
stands the multi-variable Hermite polynomials of order $|\beta|=2m$.

Hence, the hypergeometric
functions of the form ${~}_0F_1\left( s+1;-\frac{t^2}{4}\left(X_h\right)^2\right)\a$, 
provided by Corollary \ref{hypercomplexEGFCorollary}, are explicitly given by
\begin{eqnarray*}
{~}_0F_1\left( s+1;-\frac{t^2}{4}\left(X_h\right)^2\right)\a
=\sum_{m=0}^\infty \frac{\left(h^2t\right)^{2m}}{(s+1)_m} \sum_{\|x\|_1=2mh}  
H_{\frac{|x|_1}{h}}\left(-\frac{1}{2h}\e\right)\a.
\end{eqnarray*}
\end{remark}

\section{Further directions}\label{ConcludingSection}

The operational characterization given by Proposition \ref{RodriguesFormulaXh} enables us to compute,
from the knowledge
of the hypercomplex extension of the falling factorials considered in Example
\ref{HypercomplexFactorialPolynomials}, 
several classes of hypercomplex polynomials of discrete variable were derived in a simbolic way like e.g.
the hypercomplex extension of the Poisson-Charlier polynomials illustrated in Example 
\ref{HypercomplexPoissonCharlier}. Part of this characterization 
implies that for each $k\in\BN_0$ the hypercomplex Poisson-Charlier 
polynomial of degree $k$ underlying to the parameter values $a>0$
yields as a solution of the following differential-difference equation
\begin{eqnarray}
\label{evolutionPoissonCharlier} \left\{\begin{array}{lll} 
\partial_a \f(x,a)= -a\sum_{j=1}^n\partial_h^{+j}\f(x,a) & \mbox{for} & (x,a)\in
h\BZ^n \times \left[0,\infty\right)
\\ \ \\
\f(x,0)=\w_k(x;h) & \mbox{for} & x\in h\BZ^n.
\end{array}\right.
\end{eqnarray}
whereas for $a<0$ the hypercomplex Charlier polynomial of degree $k$ may be recovered from the mapping
transformation 
$\f(x,a)\mapsto \f(x,-a)$ on (\ref{evolutionPoissonCharlier}).

From a general perspective, the integral representation given by Corollary \ref{RodriguesFormulaXhCorollary}
enables us to determine each $\w_k(x;h;\lambda)$
as $\w_k(x;h;\lambda)=\int_{0}^\infty \f(x,s) ds$, where $\f(x,s)$ is a solution of the
following differential-difference time-evolution problem on the space-time domain 
$h\BZ^n \times \left[0,\infty\right)$ 
\begin{eqnarray}
\label{evolutionSigmaDh} \left\{\begin{array}{lll} 
\partial_s \f(x,s)= -s\sigma(D_h^+)\f(x,s) & \mbox{for} & (x,s)\in
h\BZ^n \times \left(0,\infty\right)
\\ \ \\
\f(x,0)=\w_k(x;h) & \mbox{for} & x\in h\BZ^n.
\end{array}\right.
\end{eqnarray}
  
Another perspective for the integral representation of $\w_k(x;h;\lambda)$
obtained in Corollary \ref{RodriguesFormulaXhCorollary} may also be obtained through the 
representation of the exponentiation operator $\exp\left(-s\sigma(D_h^+)\right)$ on
the momentum space on the $n-$cube $Q_h=\left( -\frac{\pi}{h},\frac{\pi}{h}\right)^n$, the corresponding
{\it Brioullin zone}
\footnote{Roughly speaking, the {\it Brioullin zone} 
corresponds to the cellular decomposition of $\BR^n$ determined
from a lattice such that there is a one-to-one correspondence between each point of the lattice
and each boundary point of the underlying cell on the momentum space.
}
of the lattice $h\BZ^n$. Such representation involves 
the 'discrete Fourier' \footnote{The 'discrete Fourier' transform is a particular case of a
{\it Fourier quadrature rule}
used to represent a lattice function on the momentum space. 
  Further details may concerning this may be found in \cite{Froyen89}.} transform  
  (cf.~\cite[Subsection 5.2]{GuerlebeckSproessig97} and \cite{GuerlbeckHommel03})
\begin{eqnarray*}
(\mathcal{F}_h \g)(y)&=&\left\{\begin{array}{lll} 
 \frac{h^n}{\left(2\pi\right)^{\frac{n}{2}}}\sum_{x\in h\BZ^n}\g(x)\exp(ix \cdot y) & \mbox{for} & y\in Q_h
\\ \ \\
0 & \mbox{for} & y\in \BR^n \setminus Q_h.
\end{array}\right.
\end{eqnarray*}

Here we recall that the 'discrete Fourier' transform $\mathcal{F}_h$ is a unitary operator from 
$\ell_2(h\BZ^n)$ onto $L_2(Q_h)$, whose inverse is given by the restriction of the Fourier transform 
$\mathcal{F}:L_2(\BR^n) \rightarrow L_2(\BR^n)$ to the lattice $h\BZ^n$, that is 
$\mathcal{F}_h^{-1}=\mathcal{R}_h\mathcal{F}$ where $\mathcal{R}_h \g(x)$ stands
the restriction of the function $\g(x)$
to $h\BZ^n$ and
\begin{eqnarray}
\label{FourierTransform}
(\mathcal{F} \g)(x)=\frac{1}{\left(2\pi\right)^{\frac{n}{2}}}\int_{\BR^n} \g(y)\exp(-i x \cdot y) dy.
\end{eqnarray}

On the other hand, from the unitary one-to-one correspondence $\kappa(\partial_{x_j}) \mapsto \kappa(-iy_j)$
provided by Fourier transform (\ref{FourierTransform}), one infers from the construction of $\lambda(D_h^+)$ 
and $\sigma(D_h^+)$ provided by Proposition \ref{FourierDualProp} and Proposition \ref{RodriguesFormulaXh},
respectively, the following one-to-one correspondence
$$\mathcal{F}:\exp\left(-s\sigma(D_h^+)\right)\g(x)
\mapsto \exp\left(- s\lambda\left(\frac{D_{h}^+ \exp(-i x \cdot y)}{\exp(-i x \cdot y)}\right)\right)
(\mathcal{F}\g)(y).
$$

Based on this, one can thus infer that the integral description of $\w_k(x;h;\lambda)$, 
provided by Corollary \ref{RodriguesFormulaXhCorollary}
corresponds in the momentum space to the following integral representation over $Q_h \times [0,\infty)$:
\begin{eqnarray*}
\frac{\lambda({\bf 0})}{\left(2\pi\right)^{\frac{n}{2}}}\int_{0}^\infty \int_{Q_h} 
\exp\left(- s\lambda\left(\frac{D_{h}^+ \exp(-i x \cdot y)}{\exp(-i x \cdot y)}\right)-i x\cdot y
\right) (\mathcal{F}_h\w_k)(y;h)dy ds.
\end{eqnarray*}

Integral representations of the above type may be applied
to study spectral problems in the context of discrete quantum mechanics
 (cf. \cite[Section 6]{DHS96}). 
 Although the eigenfunctions of
any finite difference Dirac operator may be computed in terms of formal power series representations 
(for example, the hypercomplex formulation of {\bf EGF} provided by equation (\ref{evolutionEGF})), 
the existence of 'doublers'\footnote{From the Fourier analysis side, the presence
of 'doublers' is equivalent to say that the symbol of the finite difference Dirac
operator has many zeroes inside 
the {\it Brillouin zone}.} 
inside the {\it Brillouin zone}
is notoriously the major difficult in case that self-adjoint discretizations such as 
$\frac{1}{2}\left( D_{h/2}^+ + D_{h/2}^-\right)$ 
are considered (cf.~\cite[Chapter 4]{Rothe05}).

Besides the hypercomplex extension of multi-variable falling factorials $(x;h)_\alpha$ 
defined in (\ref{factorialPolynomials}), the construction obtained
in Example \ref{HypercomplexFactorialPolynomials} 
gives us some interesting insights concerning applications in sampling theory and integral transforms.
Part of the construction provided on this example even shows that the hypercomplex polynomials $\w_k(x;h)$
may be reconstructed from the sampling points 
determined from the intersection of the lattice $h\BZ^n$ with the level curves determined from the 1-norm constraint 
\begin{eqnarray*}
\label{LevelCurve}\sum_{j=1}^n |x_j|=2 \left\lfloor \frac{k}{2}\right\rfloor h,
\end{eqnarray*}
as it is despicted in Figure 1.

Based on this observation, one can thus obtain the Taylor series approximation
for a continuous function $\f(x)$ defined 
on a bounded domain $\Omega$ of $\BR^n$ by simply take the $m-$term truncation of
the hypercomplex {\bf EGF} $\G_h(x,t)$
(cf.~Corollary \ref{hypercomplexEGFCorollary} and Remark \ref{hypercomplexEGFRemark}). 
The choice of $m$ is based on the fact
that the level curve labeled by the parameter $k=2m$ shall corresponds to the maximal level curve 
contained in $\Omega$. 

\begin{figure}
\centering
\includegraphics[scale=0.25]{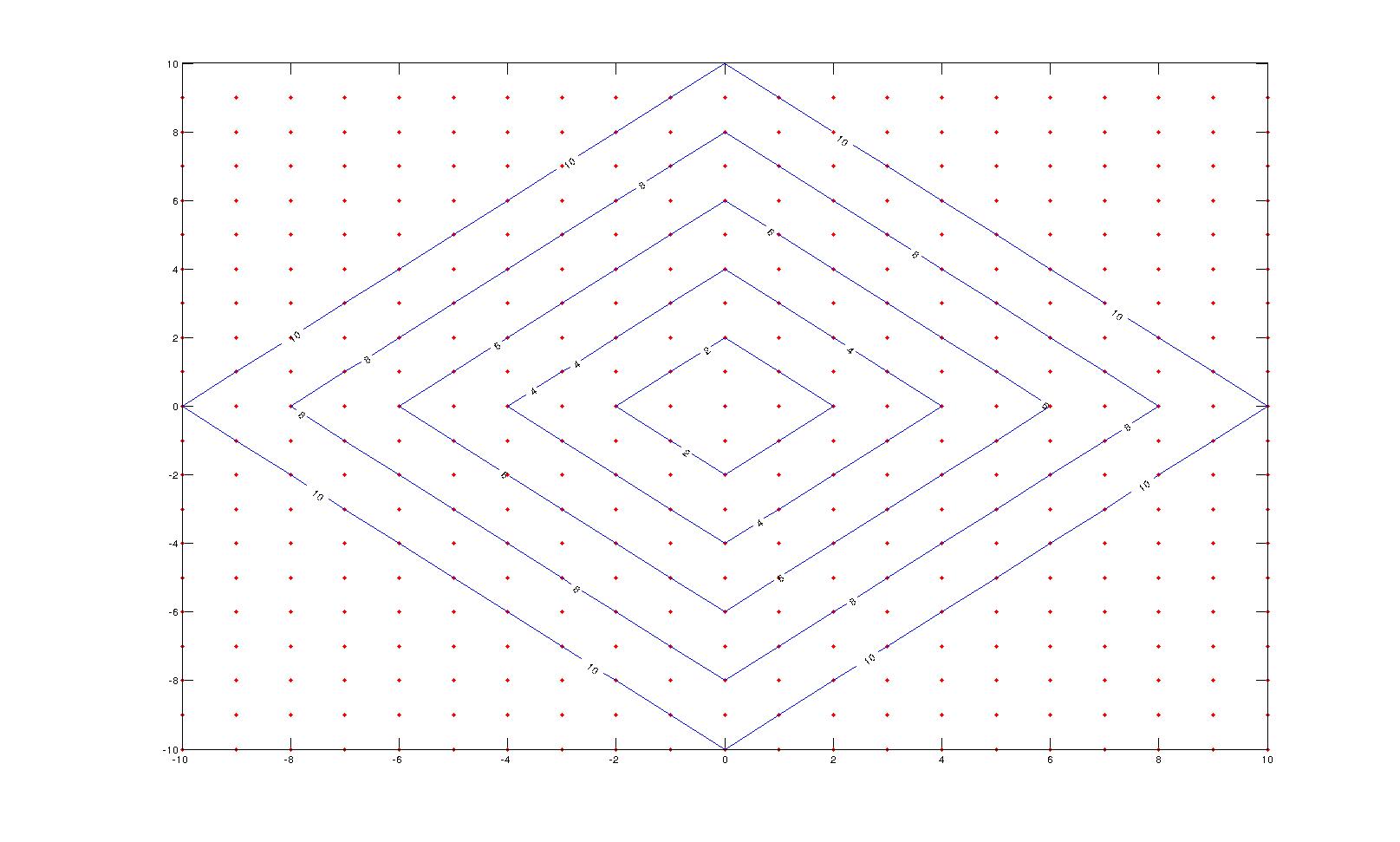}
\caption{Construction of the hypercomplex polynomials on the lattice $[-10,10]^2 \cap \BZ^2$
provided by Example \ref{HypercomplexFactorialPolynomials}.}
\end{figure}

This point of view, roughly considered in \cite{RSS13}, was applied recently in
\cite{BaaskeBernsteinDeRidderSommen14} to compute the discrete counterpart of the {\it heat kernel}
\footnote{i.e. the kernel of 
$t^{-\frac{n}{2}}(W\f)\left(\frac{x}{\sqrt{t}}\right)$, where $W$ denotes the Weierstra\ss~transform
considered in Example \ref{HypercomplexFactorialPolynomials}.} as a solution
of the differential-difference equation
\begin{eqnarray*}
\label{discreteHeatEq} \left\{\begin{array}{lll} 
\partial_t \g(x,t)= \sum_{j=1}^n \dfrac{\g(x+h\e_j,t)+\g(x-h\e_j,t)-2\g(x,t)}{h^2} &, (x,t)\in
h\BZ^n \times \left[0,\infty\right)
\\ \ \\
\g(x,0)=\m_\alpha(x) &,  x\in h\BZ^n.
\end{array}\right.
\end{eqnarray*}

On this direction one would like to stress that the construction considered in \cite{RSS13}
shall be understood as special case 
of a more general framework on which the interpolating functions 
are computed upon the action of an integral transform on the quasi-monomials $\m_\alpha(x)$ constructed from
the operational rule (\ref{quasiMonomialQuantumFieldLemma}). 
On the other hand, the action of the propagator 
$\exp\left(-\frac{t}{4}\left(D_{h/2}^++D_{h/2}^- \right)^2\right)$ on each $\m_\alpha(x)$
results into solutions of the above differential-difference equation.
Hereby each $\m_\alpha(x)$ may be taken from the range of multi-variable polynomials explored in 
Example \ref{centralDifferencesExample} (or even in \cite[Example 4 of Section 2]{DHS96})
through the substitution $h \rightarrow \frac{h}{2}$.

Obviously, this approach can be applied and extend to other shapes. 
For example, the inspiring work of Kisil \cite{Kisil02} sheds some insights about how such
{\it quasi-monomiality} formulation shall be applied to formulate path integrals. Its feasibility will be discussed
on further research in comparison with other well-known path integral formalisms (cf.~\cite[Chapter 2]{Rothe05}).

\subsection*{Acknowledgments}
The author acknowledges the anonymous referees for the various helpful comments and suggestions.


\end{document}